\documentstyle[12pt,epsf]{article}

\title{\Large \bf  Stationarity  of Switching VAR and Other Related Models}
\author{
Gopal K. Basak \thanks{Dept. of Mathematics, University of Bristol,
Bristol BS8 1TW, U.K., magkb@bris.ac.uk, and Stat-Math Unit, Indian Statistical Institute, Kolkata - 700 108, India, gkb@isical.ac.in}
\and 
Zhan-Qian Lu \thanks{
Statistical Engineering Div, ITL,
NIST,
Gaithersburg, MD 20899-8980,
U.S.A.
}
}

\renewcommand{\baselinestretch}{1.5}
\makeatletter
\@addtoreset{equation}{section}

\def\singlespace{\def\baselinestretch{1}\@normalsize}

\makeatother
\makeatletter
\renewcommand{\section}{\@startsection{section}{1}{0mm}%
{-\baselineskip}{0.3\baselineskip}{ \bf}}
\renewcommand{\subsection}{\@startsection{subsection}{1}{0mm}%
{-\baselineskip}{0.05\baselineskip}{ \bf }}
\makeatother
\date{}
\begin{document}
\maketitle

\newtheorem{theorem}{Theorem}
\newtheorem{lemma}{Lemma}
\newtheorem{corollary}{Corollary}
\newtheorem{definition}{Definition}
\newtheorem{proposition}{Proposition}
\newtheorem{conjecture}{Conjecture}

\newcommand{\be}{\begin{equation}}
\newcommand{\ee}{\end{equation}}
\newcommand{\beqa}{\begin{eqnarray*}}
\newcommand{\eeqa}{\end{eqnarray*}}
\newcommand{\beqan}{\begin{eqnarray}}
\newcommand{\eeqan}{\end{eqnarray}}
\newcommand{\bd}{\begin{description}}
\newcommand{\ed}{\end{description}}
\newcommand{\arrow}{\rightarrow}
\newcommand{\ov}{\stackrel}
\newcommand{\defi}{\stackrel{\triangle}{=}}
\newcommand{\indistr}{\stackrel{d}{=}}
\newcommand{\cind}{\stackrel{d}{\arrow}}
\newcommand{\cinp}{\stackrel{p}{\arrow}}
\newcommand{\e}{\varepsilon}
\newcommand{\norm}{{\parallel}}

\newcommand{\E}{{\rm E}}
\newcommand{\pr}{{\rm P}}
\newcommand{\smallo}{{o}}
\newcommand{\bigo}{{O}}
\newcommand{\cov}{{\rm Cov}}
\newcommand{\corr}{{\rm Corr}}
\newcommand{\var}{{\rm Var}}
\newcommand{\avar}{{\rm AVar}}
\newcommand{\acov}{{\rm ACov}}
\newcommand{\vech}{{\rm vech}}
\newcommand{\vecv}{{\rm vec}}
\newcommand{\tr}{{\rm Tr}}
\newcommand{\dh}{\hat{D}}
\newcommand{\tilt}{\tilde{T}}
\newcommand{\tmh}{\hat{T}^m}
\newcommand{\tmht}{\hat{T}^{m \prime}}

\newcommand{\bx}{{\bf x}}
\newcommand{\bu}{{\bf u}}
\newcommand{\bv}{{\bf v}}
\newcommand{\bw}{{\bf w}}
\newcommand{\bbx}{{\bf X}}
\newcommand{\by}{{\bf y}}
\newcommand{\ba}{{\bf a}}
\newcommand{\nn}{\nonumber}

\newcommand{\mt}{^{m \prime}}
\newcommand{\diag}{{\rm diag}}
\newcommand{\dg}{{\rm dg}}
\newcommand{\dgv}{{\rm dgv}}
\newcommand{\R}{{\Re}}
\newcommand{\Z}{{\cal Z}}
\newcommand{\F}{{\cal F}} 
\newcommand{\St}{{\cal S}}

\newcommand{\Th}{\Theta}
\newcommand{\th}{\theta}
\newcommand{\no}{\nonumber}
\newcommand{\La}{\Lambda}

\begin{abstract}
Switching ARMA models greatly enhance the standard linear models to 
the extent that  different ARMA model is allowed in a different 
regime, and the regime switching   is  typically   assumed 
a Markov chain on the finite states of potential  regimes.  
Although statistical issues have been the subject of many recent papers,  
there is few systematic study of the probabilistic aspects of 
this new class of  nonlinear models.  
 This paper discusses some basic issues concerning this class of 
models including strict stationarity, influence of initial conditions, 
and second-order property by studying SVAR models.
 A number of examples are given to 
illustrate the theory and the variety of applications.  
Extensions to other models such as mean-shifting, and inhomogeneous 
transition probabilities are discussed. 
\end{abstract}

\begin{singlespace}
\begin{footnotetext}
{
{\em Abbreviated title}:  Switching ARMA Models.
 \\
\hspace*{.2in} {\em AMS 1991  subject classification.} 
Primary 62M10; secondary 60G10. \\
\hspace*{.2in} {\em Key words and phrases}. 
Markov switching, strict stationarity, Lyapunov exponents, 
influence of initial conditions, mixed process, mean-shifting model. }
\end{footnotetext}
\end{singlespace}
\newpage

\section{Introduction}
\label{sec:intro}
Switching ARMA models belong to a new class of time series models 
which  are capable of capturing various nonlinear aspects of time series 
data such as nonnormality, 
asymmetry, irreversibility, and variable predictability 
[e.g. Hamilton 1989; Huges and Guttorp 1994; Krolzig 1997; 
  Lu and Berliner 1997]. 
This class of models extends the ARMA linear system to 
the extent that  different ARMA model is allowed in a different 
regime, and the regime switching   is  typically   assumed 
a Markov chain on the finite states of potential  regimes. 
While  statistical aspects of fitting these models have been 
 much  discussed as summarized by   Krolzig (1997);  
There is, however,  few systematic study of the probabilistic 
aspects of switching ARMA models, such as 
stationarity or ergodicity.

 This paper discusses some general 
conditions that  ensure stationarity and  other probabilistic 
properties such as  existence of moments. A general theory due to 
Brandt (1986) is reviewed (Section~\ref{sec:brandt}).  
A theory of stability (or, of the noninfluence of initial conditions) of
switching vector autoregressive models (SVAR) is developed 
(Section~\ref{sec:stability}).  Some interesting 
examples are given to illustrate the subtle generality of 
the developed stationarity 
conditions and the variety of applications of the switching vector 
autoregressive models. For example, 
we exhibit (as in  Holst et al (1994)) that unstable subprocesses and stable 
processes can be mixed 
to produce a stationary  process (Example 2), two unstable subprocesses 
can still be mixed to be stationary (Example 4), and stable subprocesses 
may not always produce stationary mixed process, and a counter-example 
is given (Example 3). 
The second-order theory of switching AR models is developed  
(Section~\ref{sec:second-order}). 
We also discuss the mean shifting 
models (Section~\ref{sec:shift}), 
 switching moving average, and switching ARMA 
models (Section~\ref{sec:sarma}). 

\section{General theory}
\label{sec:strict-stat}

\subsection{Switching vector AR models}
\label{sec:sar}
A general model is the following vector stochastic difference equation  
\be
X_n=A_n X_{n-1}+E_n, \mbox{  $n\in \Z$}, 
\label{eq:stochastic-eq}
\ee
where $X_n\in \R^p$ and $A_n$ is a $p\times p$ matrix and $E_n\in \R^p$ 
is a noise vector.  Various additional structure will be imposed  
 on $A_n, E_n$ later.   For example, an AR(p) process can be represented as  
(\ref{eq:stochastic-eq}) in which $A_n$ is a constant matrix assuming a 
special structure. When $\{(A_n, E_n)\}$ is iid, (\ref{eq:stochastic-eq}) 
is called the {\em Random Coefficient Autoregressive (RCA)} model 
(Nicholls and Quinn, 1992). 
Since in large part such a system  is used for modelling 
 stationary time series data,  stationarity property is 
a priority in the study of probabilistic aspects of such 
random dynamical systems. 
A  theory for the general stochastic equations (\ref{eq:stochastic-eq}) 
 is reviewed in Section~\ref{sec:brandt}. 
However, one of our objectives
 is to study the  so-called 
{\em Markov switching vector AR(1) model} (SVAR(1)):  
Suppose there are $r$
 potential regimes, say $\St=\{1, 2, \ldots, r\}$ and 
$I_n$ is a Markov chain taking on values in $\St$. 
Define 
\be
A_n=\sum_{i=1}^r B_i 1_{\{I_n=i\}},
\label{eq:coefficient in svar}
\ee
where $B_1,\ldots, B_r$ are $r$ unknown  or partially unknown 
$p\times p$ matrices; and  
\be
E_n=\sum_{i=1}^r \Sigma_i \e_{ni} 1_{\{I_n=i\}}
\label{eq:noise in svar}
\ee
where $\{ \e_{ni}\}$ are independent processes,  each 
subsequence is iid within itself,  having zero mean  and 
identity covariance  matrix. 
In addition, we make the assumption of {\em independence},
 that $\{I_n\}$ is independent of 
noise processes $\{\e_{n1}, \e_{n2},\ldots, \e_{nr}\}$.
We also assume that $\{I_n\}$ is {\em irreducible} and {\em aperiodic}, 
thus {\em ergodic}.

First we ask the question whether there exists a 
strict stationary solution for (\ref{eq:stochastic-eq})? 
Since a strict stationary process may not have  any moment existing, 
this is a fairly weak assumption. 
Though  necessary and sufficient stationarity conditions for 
 RCA models  are available  [e.g. Nicholls and Quinn (1992)], 
necessary and sufficient stationarity conditions when $\{A_n\}$ is dependent 
has not yet been given (see, however, Bougerol and Picard 1992). 
However, very general sufficient condition
that ensure stationarity  can be formulated from the proof of  
Brandt (1986), and is first made known in  Bougerol and Picard (1992).
For convenience of later use, we will  
restate a general theorem related to this theory.  
Here it will be assumed  that the super-process 
$\{(A_n, E_n)\}_{n=-\infty}^{\infty}$ 
are (jointly) stationary matrices and 
vectors. It appears that all known results in this area 
make this convenient assumption, though more can be said 
in our setup (later). 

\subsection{Brandt's 
result
}
\label{sec:brandt}
We first state a general 
result
 giving sufficient conditions 
for strict stationarity. 
Here, we do {\em  not} need to assume that $A_n$ takes on 
discrete values  $B_i$'s.  In the case of SVAR(1), 
stationarity of $\{A_n, E_n\}$ is equivalent to  
assuming that  the ergodic  chain $\{I_n\}$ 
starts from the remote past or $I_0$ takes on the stationary 
distribution.

The tool is the theory of {\em Lyapunov exponents} or 
product of random matrices. A technical assumption 
that ensures  existence of Lyapunov exponents 
for a stationary   sequence of  random matrices  
$A_1, A_2,\ldots, A_n,\ldots$ is 
\be
\E \max(\log \|A_1\|,0)<\infty.
\label{eq:le-cond}
\ee
This is obviously satisfied if $A_1$ takes on only finite number 
of values as in the case of SVAR(1). 

Under (\ref{eq:le-cond}) the (largest) Lyapunov exponent  
is defined as  
\be
\lambda=\lim_{n\arrow \infty} (1/n)\log \|A_n\ldots A_1\|
\label{eq:lyapunov-exp}
\ee
which holds almost surely.

Furthermore, if the process is ergodic, 
the Lyapunov exponent is constant and 
\be
\lambda=\inf \{(1/n) \E\log\|A_n\ldots A_1\|, n\geq 1\}.
\label{eq:lambda}
\ee
The existence of the limit in  
(\ref{eq:lyapunov-exp}) and (\ref{eq:lambda}) can be justified  using    
Kingman's subadditive ergodic theorem.
Similarly, the following limit theorem 
holds with the same Lyapunov exponent as a consequence of 
(\ref{eq:lambda}):
\be
\lambda=\lim_{n\arrow \infty} (1/n)\log \|A_0 A_{-1}\ldots A_{-n+1}\|.
\label{eq:lyapunov-exp2}
\ee
Note that $\lambda$ is defined {\em independent} of the particular 
matrix norm used.

More generally, under stationarity and (\ref{eq:le-cond}), 
one can apply  Oseledec's multiplicative ergodic theorem 
to define a spectrum of Lyapunov exponents 
$\lambda_1=\lambda\geq \lambda_2\geq \ldots \geq \lambda_p$:  
\be
\lambda_i=\lim_{n\arrow \infty} \frac{1}{n} \log \delta_i(n), 
\mbox{ holds almost surely  for $1\leq i\leq p$},    
\ee
where $\delta_1(n)\geq \ldots \geq \delta_p(n)$ 
are the singular values of $A_n A_{n-1} \ldots A_1$.
Under ergodicity, $\lambda_i$'s are constants, independent of the 
particular realization in $\{A_n\}$.

\begin{proposition} \label{th:strict-stat}
Given that the super-process $\{A_n,E_n\}$ is stationary and ergodic. 
Suppose that $P(A_0=0)>0$ or the following conditions are met:  
(\ref{eq:le-cond}) holds and the Lyapunov exponent for $\{A_n\}$  is 
negative; that is
\be
\mbox{(NL): \hspace{.9in} } 
\lambda=\lim_{i\arrow \infty}^{} (1/i) \log\|A_0A_{-1}\ldots A_{-i+1}\|
<0
\label{eq:lyapunov-neg}
\ee
 and the noise  satisfies 
\be
\mbox{$\E \max(\log\|E_{1}\|,0)<\infty$.} 
\label{eq:noise-cond}
\ee
Then (i) 
\be
W_n=E_n+\sum_{i=0}^{\infty} A_n A_{n-1} \ldots A_{n-i} E_{n-i-1}
\label{eq:limiting-proc}
\ee
 is the only proper stationary solution of 
(\ref{eq:stochastic-eq}) for the given $\{A_n,E_n\}$. 
(ii) The sum on the right-hand side of (\ref{eq:limiting-proc})  
converges absolutely almost surely. 
(iii) 
Furthermore, 
\be
P(\lim_{n\arrow + \infty} | X_n(\bx)- W_n|=0)=1,
\ee
for arbitrary random variable $X_{-m-1}=\bx$ at time $-m-1$ 
(defined on the same probability 
space as $\{A_n, E_n\}$), in particular
\be
X_n(\bx)\cind W_0, \mbox{ as $n \arrow + \infty$.}  
\ee
\end{proposition}

Part (i) of this 
result
is first given in Bougerol and Picard (1992). 
The proof is similar to the one-dimensional case as proved by Brandt (1986) 
under a stronger assumption. See Bougerol and Picard (1992) for 
more details and a necessary condition. 

For convenience of readers, we prove the following lemma for part (ii). 
\begin{lemma} \label{la:converg}
If the stationary super-process $\{A_n, E_n\}$ satisfies 
(\ref{eq:lyapunov-neg}) and (\ref{eq:noise-cond}), 
 the RHS of (\ref{eq:limiting-proc}) 
converges  absolutely  almost surely. 
\end{lemma}
{\em Proof.} 
First, by (\ref{eq:lyapunov-neg}) and (\ref{eq:noise-cond})   
\beqa
\lefteqn{
 \limsup_{i\arrow \infty} 
 (1/(i+1)) \log \|A_n A_{n-1} \ldots A_{n-i} E_{n-i-1}\|}\\
&\leq& \limsup_{i\arrow \infty} (1/(i+1))\log \|A_n A_{n-1} \ldots A_{n-i}\|+
(1/(i+1))\log \|E_{n-i-1}\| \\
&\leq & \lambda+0 <0, \mbox{ a.s.}
\eeqa
which implies that 
\[ \limsup_{i\arrow \infty} 
 \|A_n A_{n-1} \ldots A_{n-i} E_{n-i-1}\|^{1/i}<1\mbox{ a.s.}.\]
Thus, 
 the RHS of (\ref{eq:limiting-proc}), which is bounded by 
\[ \|E_n\|+\sum_{i=0}^{\infty} 
\|A_n A_{n-1} \ldots A_{n-i} E_{n-i-1}\|,
\]
is absolutely convergent almost surely by virtue of 
Cauchy's root criterion.
\hfill $\Box$

Since the process $W_n$ defined by 
(\ref{eq:limiting-proc}) is a well-defined  moving average  
function of ergodic stationary process $\{A_n, E_n\}$, 
it follows that it is stationary and ergodic. 
Thus, $W_n$ is a $MA(\infty)$ process with random coefficients.

A key idea in the proof of Proposition~\ref{th:strict-stat} is based on 
 the following expansion  which holds  for any integers $m$ and $n$ 
as implied by the recursive nature of (\ref{eq:stochastic-eq})  
\beqan
X_n(\bx) &=& A_n A_{n-1}\ldots A_1A_0A_{-1}\ldots A_{-m} \bx 
 \nonumber\\
    &&+\sum_{i=0}^{n+m-1} A_nA_{n-1}\ldots A_{n-i}E_{n-i-1} +E_n
\label{eq:expansion}
\eeqan 
where $X_n(\bx)$ can be interpreted as the state at time $n$  of 
the system governed by (\ref{eq:stochastic-eq})  
if it starts at time $-m-1$ with the random 
initial state $X_{-m-1}=\bx$.
Thus, $W_n$ can be regarded as the limit of 
$X_n(\bx)$ starting from the remote past. 

Further, part (iii) of the theorem says that $X_n(\bx)$ converges 
to $W_n$ forward in time as time $n$ tends to the future. 
This follows from that 
\beqa
X_n(\bx)-W_n &=&A_n A_{n-1}\ldots A_1A_0A_{-1}\ldots A_{-m} \bx\\
             && + \sum_{i=n+m}^{\infty} A_n A_{n-1}\ldots A_{n-i} E_{n-i-1}
\eeqa
which tends to zero almost surely under  condition (NL) 
thanks to   Lemma~\ref{la:converg}. 

{\em Remark} 1. Since for any positive random variable $X$, 
by  Jensen's  inequality   $E \log X\leq \log E X$ holds 
whenever $E X<\infty$, it follows that  
 {\em whenever  $E X^{\alpha}<\infty$ for any $\alpha>0$} 
we have 
$E \log X<\infty$ and hence $E \max(0,\log X)<\infty$. 
(Note that $\max(0,\log X)$ represents the positive 
part of $\log X$.)

Next, we consider the more realistic situation that a Markov 
switching process starts from a {\em finite} time in the past and 
discuss when such a process can be stationary and ergodic. 

\subsection{Stability of SVAR models}
\label{sec:stability}  

Under (\ref{eq:le-cond}) the (largest) Lyapunov exponent  
is defined as in  \ref{eq:lyapunov-exp}
\be
\lambda=\lim_{n\arrow \infty} (1/n)\log \|A_n\ldots A_1\| \ \ \mbox{almost surely.}
\ee

Now consider the situation that 
the SVAR(1) process starts at some {\em fixed}  time, say time $0$,  
with some {\em arbitrary} starting value $X_0$  and the regime 
process $\{I_n\}$ starts from an {\em arbitrary} distribution $I_0$. 
Let $X_n(X_0, I_n(I_0))$ denote the process evolved according to 
(\ref{eq:stochastic-eq}) with starting value $X_0$   
and starting regime $I_0$ at time $0$.  
The question arises as to what's the influence 
of the initial condition or the transient effect. Naturally, 
one would hope that the initial effect will eventually be washed out 
or vanish. 
It is indeed so. We prove it in the next theorem
after illustrating a lemma.
The result of this lemma is well known (e.g., Bhattacharya and Waymire (1990),
 p.197) however we put it here for the sake of completeness.

\vspace{.5em}
\noindent
\begin{lemma} \label{tm:exp}
{\it Let $I^{1} (\cdot)$ and $I^{2} (\cdot)$ be two independent 
replicas
of an irreducible and aperiodic
 Markov chain $I(\cdot)$ with finite state space $S$
(with $r$ number of elements),
having the same transition probability $(\ P = ((p_{ij}))\ )$.
Define,
$$
\tau = \inf \{ k \ge 0 \ : \ I^{1}_{k} = I^{2}_{k} \} \ .
$$
Then, for any $i, j \in S$,
 $\ P(\tau > n \ | \ I^{1}_{0} = i, I^{2}_{0} = j )$ converges to zero,
exponentially fast, as $t \to \infty$.}
\end{lemma}

\noindent
{\sc Proof.}
 Define,
$$
p(r_0) = \max_{k,l} P(I^{1}_{m} \neq I^{2}_{m},
 1 \le m \le r_0,
 \ | \ I^{1}_{0} = k, I^{2}_{0} = l ) 
$$
Since the state space is finite, under
 the condition of irreducibility and aperiodicity,
it is clear that, there exists an $r_0 \ge 1$ such that $p^{(r_0)}_{ij} > 0$ 
for all $i, j \in S$. Let $\alpha_0 = \min_{i,j} p^{(r_0)}_{ij}$. 
Then $\alpha_0 > 0$ and 
$p(r_0) \le \max_{k,l} P(I^{1}_{r_0} \neq I^{2}_{r_0},
 \ | \ I^{1}_{0} = k, I^{2}_{0} = l )
 = \max_{k,l} (1 - \sum_{i} P(I^{1}_{1} = i  \ | \ I^{1}_{0} = k))
P(I^{2}_{1} = i \ | \ I^{2}_{0} = l )
= \max_{k,l} (1 - \sum_{i} p_{ki} p_{li}) 
\le (1 - r\alpha_0^2) 
< 1
$
 Let $n \ge 1$ be an integer.
 Then, using
Markov property and stationarity of the joint Markov chain 
$(I^{1}, I^{2})$ we obtain,

$$\begin{array}{lll}
& & P(I^{1}_{m} \neq I^{2}_{m}, 1 \le m \leq nr_0
 \ | \ I^{1}_{0} = i, I^{2}_{0} = j ) \no\\
& = & 
\displaystyle \sum_{k \neq l} P(I^{1}_{m} \neq I^{2}_{m},
 1 \le m < (n - 1)r_0, 
 I^{1}_{(n - 1)r_0} = k, I^{2}_{(n - 1)r_0} = l
 \ | \ I^{1}_{0} = i, I^{2}_{0} = j ) \no\\
& & \times P(I^{1}_{m} \neq I^{2}_{m}, (n - 1)r_0 < m \leq nr_0
 \ | \ I^{1}_{(n - 1)r_0}  = k, I^{2}_{(n - 1)r_0} = l ) \no\\
& = &
\displaystyle \sum_{k \neq l} P(I^{1}_{m} \neq I^{2}_{m},
 1 \le m < (n - 1)r_0, 
I^{1}_{(n - 1)r_0} = k, I^{2}_{(n - 1)r_0} = l
 \ | \ I^{1}_{0} = i, I^{2}_{0} = j ) \no\\
& & \times P(I^{1}_{m} \neq I^{2}_{m}, 1 \le m \leq r_0
 \ | \ I^{1}_{0}  =  k, I^{2}_{0} = l ) \no\\
& \leq &
 \displaystyle \sum_{k \neq l} P(I^{1}_{m} \neq I^{2}_{m},
 1 \le m < (n - 1)r_0, 
 I^{1}_{(n - 1)r_0} = k, I^{2}_{(n - 1)r_0} = l
\no\\
&&
 \ | \ I^{1}_{0} = i, I^{2}_{0} = j ) \times p(r_0)  \no\\
& = & 
P(I^{1}_{m} \neq I^{2}_{m},  1 \le m \leq (n - 1)r_0 
\ | \ I^{1}_{0} = i, I^{2}_{0} = j ) \times p(r_0) \ .
\end{array}
$$
Using the above argument recursively we get
$$
 P(I^{1}_{m} \neq I^{2}_{m}, 1 \le m \leq nr_0
 \ | \ I^{1}_{0} = i, I^{2}_{0} = j ) \leq p^{n}(r_0) \ .
$$
Consequently, we obtain, for any $n \ge r_0$,

\vspace{.5em}

\noindent$ P( \tau > n \ | \ I^{1}_{0} = i, I^{2}_{0} = j )$
$$\begin{array}{lll}
&=& P(I^{1}_{m} \neq I^{2}_{m}, 1 \le m \leq n
 \ | \ I^{1}_{0} = i, I^{2}_{0} = j ) \no\\
 & \leq & 
 P(I^{1}_{m} \neq I^{2}_{m}, 1 \le m \leq [n/r_0]r_0
 \ | \ I^{1}_{0} = i, I^{2}_{0} = j ) \no\\
& \leq &
 p^{[n/r_0]}(r_0) \ ,
\end{array}
$$
where $[t] =$ the largest integer that is less than or equal to $t$.
Hence the result.\hfill{$\fbox{}$}

\begin{theorem}\label{th:stb-svar}
As in the condition (NL)
assume that 
 under (\ref{eq:le-cond}) the (largest) Lyapunov exponent $\lambda$,  
defined as,  
\be
\lambda := \lim_{n\arrow \infty} (1/n)\log \|A_n\ldots A_1\| < 0 \ \ \mbox{almost surely}.
\label{eq:lyapunov-exp_cond}
\ee

Under this assumption the SVAR process is stable, i.e., it has unique asymptotic 
distribution that is free from the influence of the initial distribution.
\end{theorem}

\noindent
{\bf Proof.}
Let us assume first that $\{I_n\}$ starts at $I_0$ which is the stationary distribution for
the ergodic chain.
Then, it follows that, $\{A_n, E_n\}$ are stationary. Hence
\beqan
X_n(X_0, I_n(I_0)) 
&=& A_n A_{n-1}\ldots A_1 X_0
+\sum_{i=0}^{n-1} A_nA_{n-1}\ldots A_{n-i}E_{n-i-1} +E_n \nonumber\\
&=& A_n A_{n-1}\ldots A_1 X_0
+\sum_{i=0}^{n-1} A_{i+1}A_{i}\ldots A_{1}E_{0} +E_0 \ \ \mbox{in distribution.} \nonumber\\
\label{eq:svarexpansion}
\eeqan 
Then
for any fixed $i \ge 0$,  
\beqan
\lefteqn{
 \limsup_{n\arrow \infty} 
 (1/(i+1) \log \|A_{i+1} A_{i} \ldots A_{i} E_{0}\|} \nonumber\\
&\leq& \limsup_{n\arrow \infty} (1/(i+1))\log \|A_{i+1} A_{i} \ldots A_{1}\|+
(1/(i+1))\log \|E_{0}\| \nonumber\\
&\leq & \lambda+0 <0, \mbox{ a.s.}
\label{eq:svarlyapounov}
\eeqan
which implies that 
\[ \limsup_{i\arrow \infty} 
 \|A_{i+1} A_{i} \ldots A_{1} E_{0}\|^{1/(i+1)}<1\mbox{ a.s.}.\]
Thus, 
 the RHS of (\ref{eq:svarexpansion}) is bounded by 
\[ \|E_0\|+\sum_{=0}^{\infty} 
\|A_{i+1} A_{i} \ldots A_{1} E_{0}\|,
\]
which is absolutely convergent almost surely by  
Cauchy's root criterion
and $\|A_{n} A_{n-1} \ldots A_{1} X_0\| \arrow 0$ as $n \arrow \infty$ for any
$X_0$ as in (\ref{eq:svarlyapounov}).
Therefore,
$X_n(X_0, I_n(I_0))$ converges in distribution as $n \arrow \infty$ whenever $I_0$
starts from the stationary distribution. 

Let us now observe,
\beqan\label{eq:svar0}
X_n(X_0, I_n(I_0)) - X_n(X_0', I_n(I_0)) 
& = & A_n (X_{n-1}(X_0, I_n(I_0)) - X_{n-1}(X_0', I_n(I_0))) \nn\\
& = & \cdots = A_n A_{n-1} \ldots A_1 (X_0 - X_0') 
\eeqan
Thus, we obtain,
\[ (1/n) \log (|X_n(X_0, I_n(I_0)) - X_n(X_0', I_n(I_0))|)\]
\[\leq 
 (1/n) \sum_{j=1}^{n} \log (\| A_j\|)  + (1/n) \log(|(X_0 - X_0')|)\]
and hence by strong law for $\{A_j\}$'s,
and under the condition 
(\ref{eq:lyapunov-exp_cond}),
we obtain 
that the 
 distance between $X_n(X_0, I_n(I_0))$ and $X_n(X_0', I_n(I_0))$
 converges to zero, almost
 surely, exponentially fast regardless of $I_0$ as $n$ tends to infinity. 

To see that  $X_n(X_0, I_n(I_0))$ and $X_n(X_0',I_n(I_0'))$ have 
same asymptotic distribution, 
it is important to notice that, for $I_n(I_0)$ and $I_n(I_0')$
two independent finite state ergodic Markov chain starting at $I_0$ and
$I_0'$ respectively, will meet at some finite stopping time, say $\tau$,
(whose all moments are also finite) with probability one.

Define,
\[
\tilde{I_n}(I_0') = \left\{ \begin{array}{llll}
                   I_n(I_0'), & & {\rm for\  } n < \tau \\
                   I_n(I_0),  & & {\rm for\  } n \geq \tau,
                           \end{array}
\right.
\]
i.e., 
$\tilde{I_n}(I_0')$  follows
the chain $I_n(I_0')$ in the beginning and switches to $I_n(I_0)$ at the
stopping time $\tau$ moves along the same path thereafter.
Since
$I_n(I_0')$ and the $\tilde{I_n}(I_0')$ have same initial distribution
and the transition law and hence they have the same
distribution.
Hence for any bounded and Lipschitzian $f$ we get,
\beqan\label{eq:svar1}
& & |E f(X_n(X_0, I_n(I_0))) - E f(X_n(X_0', I_n(I_0')))| \nn\\
& = & |E f(X_n(X_0, I_n(I_0))) - E f(X_n(X_0', \tilde{I_n}(I_0')))| \nn\\
& = & |E([f(X_n(X_0, I_n(I_0))) - f(X_n(X_0', \tilde{I_n}(I_0')))] I_{\tau \le 
m} ) \nn\\
& & + E([f(X_n(X_0, I_n(I_0))) - f(X_n(X_0', \tilde{I_n}(I_0')))] I_{\tau > m} 
)|
\nn\\
& \le &
|E[E([f(X_n(X_0, I_n(I_0))) - f(X_n(X_0', \tilde{I_n}(I_0')))] I_{\tau \le m} \ 
| \ {\cal F}_{\tau_m})] | \nn \\
&&+ 2 \| f \| P(\tau > m) ,
\eeqan
where $\tau_m = \tau\wedge m$ and $ {\cal F}_{j}$ is an appropriate filtration.
with respect to which $\{ I_ns, X_ns \} $ are adapted.
We restrict the class of $f$ such that the lipschitzian constant is bounded
by one and the $\|f\| \le 1$ and call that restricted class as $BL$.
Then by Markov property we get, for $m < n$,
\beqan\label{eq:svar2}
& & |E([f(X_n(X_0, I_n(I_0))) - f(X_n(X_0', \tilde{I_n}(I_0')))] I_{\tau \le m} 
\ | \ {\cal F}_{\tau_m})] | \nn\\
& = &
|E[f(X_{n-\tau_m}(z, I_{n-\tau_m}(J))) - f(X_{n-\tau_m}(z', I_{n-\tau_m}(J)))] | 
\nn\\
& \le &
E(|X_{n-\tau_m}(z, I_{n-\tau_m}(J)) - X_{n-\tau_m}(z', I_{n-\tau_m}(J)) | \wedge 
2 )
\eeqan
conditionally on $z = X_{\tau_m}(X_0, I_{\tau_m}(I_0))$,
$z' = X_{\tau_m}(X_0', I_{\tau_m}(I_0'))$ and 
$J = I_{\tau_m}(I_0)$.
Since, by earlier argument, for each $z, \ z', \ J$,
$|X_{n-\tau_m}(z, I_{n-\tau_m}(J)) - X_{n-\tau_m}(z', I_{n-\tau_m}(J)) |$
goes to zero almost surely, exponentially fast, as $n \to \infty$,
by Lebesgue's dominated convergence theorem
$E(|X_{n-\tau_m}(z, I_{n-\tau_m}(J)) - X_{n-\tau_m}(z', I_{n-\tau_m}(J)) | 
\wedge 2 ) \to 0$,
 as $n \to \infty$, almost surely, for each fixed $m \ge 1$.
Therefore, again using Lebesgue's dominated convergence theorem and
the fact that $\tau$ is finite with probability one
 (by Lemma \ref{tm:exp}),
 we obtain,
first by taking limit $n \to \infty$ and then $m \to \infty$,
\beqan\label{eq:svar3}
& & |E f(X_n(X_0, I_n(I_0))) - E f(X_n(X_0', I_n(I_0')))| \nn\\
& \leq & 
|E[E([f(X_n(X_0, I_n(I_0))) - f(X_n(X_0', \tilde{I_n}(I_0')))] I_{\tau \le m}\ 
| \ {\cal F}_{\tau_m})] |
+ 2 \| f \| P(\tau > m)  \nn\\
& & \to 0 ,
\eeqan
uniformly over bounded Lipschitzian $f$ in $BL$.
Since the class of BL
 characterizes the 
weak convergence,
and hence the theorem
(for an analogous result in continuous time, see
 Basak, Bisi and Ghosh (1999)). \hfill{$\fbox{}$}

\begin{corollary}
Under a useful and simpler condition  
where
the random matrix $A_1$ satisfies 
\beqan
\mbox{(CB): \hspace{1.4in}  } & \E \log\|A_1\|<0.
\label{eq:suff-cond}
\eeqan
for a given norm $\| \cdot \|$,
the SVAR process is stable, i.e., it has unique asymptotic 
distribution that is free from the influence of the initial distribution.
\end{corollary}

\bigskip
{\bf Proof.}
By definition (\ref{eq:lambda}), condition (CB) implies the 
negative Lyapunov condition (NL) in Proposition~\ref{th:strict-stat} 
for {\em any} norm.
Hence the proof.

\bigskip
{\bf Remark.}
Brandt (1986) focuses mainly on (CB). 
However, being independent of a matrix norm,  condition (NL) of 
Proposition~\ref{th:strict-stat} 
is  more natural  in multidimensional systems.  

\bigskip
{\bf Remark.}
It is clear that, if the assumption of irreducibility is dropped
 then one needs to restrict
attentions within the irreducible subclasses. Within each irreducible
subclass the above result is true under aperiodicity.
Also, it is easy to see, if the assumption of aperiodicity is dropped then
the above theorem fails, i.e., asymptotic distribution would have
the influence of initial distribution.

Importance of  Theorem~\ref{th:stb-svar} is in realizing the fact that
in practice, we don't have data that starts from $-\infty$ or
follows a nice initial distribution (such as the stationary
distribution), rather we have  data which starts from a finite time in the
past and with an arbitrary  initial distribution, 
usually unknown. In such a case, having
a common limiting distribution in forward time is a necessity 
in making inference of the data. 

Certainly, the question remains
in determining the rate
of convergence to the limiting distribution. 
A more interesting and challenging problem is to check for 
stability using the Lyapunov exponent approach. For this, 
a theoretical question arises: whether the analogue of 
 Kingman's subadditive ergodic theorem or more 
Oseledec's multiplicative ergodic theorem is true when the 
sequence of random matrices $\{A_n\}$ 
follows a Markov chain and the initial 
value is {\em arbitrary}? We think this is likely the case (recall 
the law of large numbers for Markov chain) but haven't seen 
any known result on this.

\section{Examples}
\label{sec:examples}
Proposition~\ref{th:strict-stat} gives a general criterion for 
checking stationarity of switching autoregressive models via  
negativity of the largest Lyapunov exponent. Theorem~\ref{th:stb-svar} 
proves the more relevant stability property under a stronger condition.   
Technique for calculating 
Lyapunov exponents for  a sequence of random matrices  
becomes very important in checking for stationarity. Unfortunately, it is 
 extremely difficult to have explicit formula of 
Lyapunov exponents except in very special cases, and 
in the general case we may have to resort to  numerical method.  

\subsection{Cases when $A_i$'s commute}
In the special cases when formula for Lyapunov exponents is available, 
condition for stationarity follows immediately. Some situations are discussed 
next. Let $A_1, A_2, \ldots$ be an ergodic stationary 
 sequence of $p\times p$ random matrices and denote $A_k=(a_{ij}(k))$.
\begin{lemma} \label{la:le-special1} 
(i) \hspace{3pt}
 If $A_k$'s are upper triangular, i.e. $a_{ij}=0$ for any $i>j$, 
and assume that $E\max(0,\log|a_{ii}|)<\infty$ for all $1\leq i\leq p$.
Then the  Lyapunov exponents exist, and they 
correspond to the ordered sequence of the 
$r$ quantities defined by 
\[ \lim_{n\arrow \infty} \frac{1}{n} \sum_{k=1}^n \log | a_{ii}(k)|
=E \log |a_{ii}(1)|, \mbox{ for $i=1,\ldots, p$.}\]
(ii) \hspace{3pt} If any pairs of matrices $A_k$'s commute,
let $\delta_1(1)\geq\ldots\geq \delta_p(1)$ be the ordered  
eigenvalues of $A_1$ and assume $E\max(0,\log |a_{11}|)<\infty$.
Then, the Lyapunov exponents exist and are given by $\lambda_i=E\log | 
\delta_i(1)|$ for $i=1,\ldots, p$. 
\end{lemma}

We now specialize the preceding theory to the switching AR model 
(\ref{eq:stochastic-eq}) when $A_n$ takes on one of 
the $r$ possible matrices $B_1, \ldots, B_r$. 
Obviously, if the sequence $A_1,A_2, \ldots $ is stationary, 
Lyapunov exponents always {\em exist}, because (\ref{eq:le-cond}) 
holds automatically. In particular, let the stationary 
distribution of $I_n$ be ${\boldmath \rho}$  such that 
$P(I_n=i)\arrow \rho_i$ for $1\leq i\leq r$ and $\rho_1+\ldots+\rho_r=1$. 
Let $E$ denote the expectation over the joint product space of 
$\{I_n\}$ and $\{\e_{ni},i=1,\ldots, r\}$ under ${\boldmath \rho}$. 
Then, (\ref{eq:lambda}) implies that 
\be 
\lambda\leq \sum_{i=1}^r \rho_i \log \| B_i\|.
\label{eq:le-ieq}
\ee
Thus, if there {\em exists a norm} such that $\|B_i\|\leq 1$ for 
$1\leq i\leq r$ where inequality holds for at least one $i$, then 
the negative Lyapunov condition is satisfied. If 
\be 
E\max(0,\log \|\e_{1i}\|)<\infty \mbox{ for $1\leq i\leq \infty$,}
\label{eq:error-cond}
\ee 
by Proposition~\ref{th:strict-stat}, the Markov switching  AR model 
with at most random walk type 
nonstationarity in subprocesses and at least one stable 
subprocess is stationary. By now  we have used the term stable process or 
stability in several places. What we mean  is the processes 
starting from  different initial conditions converge. 
In the case of a vector AR(1) process, 
this is equivalent to the coefficient matrix $A$ having eigenvalues 
whose norms are all less than one. And the latter coincides with the 
stationarity condition (cf. Example 1).

{\em Example 2.}
In  the one-dimensional case, 
negative Lyapunov condition reduces to 
$\E\log |a_n|<0$. In particular, if $a_n$ takes on 
finite numbers $b_1,\ldots, b_r$, this is 
\be
\sum_{i=1}^r \log |b_i| Pr(a_n=b_i)<0.
\label{eq:1-dim-cond}
\ee 
This is satisfied if one $|b_i|<1$ and all other 
$|b_j|\leq 1 (j\neq i)$. 
That is, under (\ref{eq:error-cond}) a switching autoregressive 
model is stable as long as 
it has a positive probability of being in a stable regime  while 
all other regimes are either stationary or random-walk type 
nonstationary.
Obviously, explosive behavior ($|b_i|>1$) in some regimes is also allowed 
as long as (\ref{eq:1-dim-cond}) is satisfied. 
\hfill $\Box$

The conclusion of Example 2 in the one-dimensional case, 
though benign and reasonable, cannot be extended to  multi-dimensional 
case, except in trivial cases such as  Lemma~\ref{la:le-special1} 
when $B_i$'s are either triangular or commutable.
Initially, we thought that the mixture of two stable processes is always 
stable. This turns out not to be true in the  multidimensional case.  
A counterexample (Example 3) is given to 
show that two stable subprocesses 
can be mixed to produce a unstable switching process. 
On the other hand, two unstable subprocesses can be 
mixed to produce a stable switching process (Example 4).

\subsection{Calculating Lyapunov exponents in a nontrivial case} 
For Example 3,  we need a result on an explicit formula 
for Lyapunov exponent in a nontrivial case due to 
Pincus (1985), see  Lima and Rahibe (1994). 
Consider the case $r=2$ and two $2\times 2$ real matrices $B_1$ and $B_2$, 
where $B_1$ is singular. 
Denote the transition probability matrix of $\{I_n\}$ by 
$\pr(I_n=j| I_{n-1}=i)=p_{ij}, i,j=1,2$ and initial 
distribution $\pr(I_0=i)=p_i,i=1,2$.

By change of basis, we can assume 
that $B_1$ takes on the form 
\[B_1=\left( \begin{array}{cc} \delta & 0\\ 0 & 0 \end{array} 
    \right)
\]
(Another form of $A$ 
\[B_1=\left( \begin{array}{cc} 0 & \delta \\ 0 & 0 \end{array} 
    \right)
\]
is not interesting  because $B_1^2=0$.)

We write $B_2^n$ in the form 
\[B_2^n=\left( \begin{array}{cc} b_{11}(n) & b_{12}(n) \\
                                 b_{21}(n) & b_{22}(n) \end{array} \right)
\]
then a result due to Pincus (1985) and  Lima and Rahibe (1994) says that 
the Lyapunov exponent is given by 
\be
\lambda=\frac{p_{21}}{p_{21}+p_{12}} \log |\delta|+
\sum_{i=1}^{\infty} p_1 p_{21} p_{12} p_{22}^{i-1} \log |b_{11}(n)|.
\label{eq:lambda-form1} 
\ee
In the case that $B_2$ is singular,  
we  consider the case that 
\[B_2=Q^{-1} \left( \begin{array}{cc} \delta_2 & 0 \\ 0 &  0
                    \end{array} \right) Q
\]
where $Q$ is an invertible matrix. 
(By a simple argument, in the other case 
$B_2=Q^{-1} \left( \begin{array}{cc} 0 &\delta_2  \\ 0 &  0
                    \end{array} \right) Q$, 
we have $\lambda=-\infty$. Not what we want.)

Then, from Lima and Rahibe (3.2), 
\be
\lambda=\frac{p_{21}}{p_{21}+p_{12}} \log |\delta| 
+\frac{p_{12}}{ p_{21}+p_{12}} \log |\delta_2|
+\frac{p_{12} p_{21}}{p_{12}+p_{21}} \log |\frac{b_{11}}{\tr(B_2)}|.
\label{eq:lya-form2}
\ee

{\em Example  3.}  Consider $B_1=\left( \begin{array}{cc} \delta_1 & 0 
                     \\ 0 & 0 \end{array} \right)$ and 
                $B_2=\left( \begin{array}{cc} b_1 & -c b_1 \\
                       b_2 & -c b_2 \end{array} \right)$.
  The eigenvalues for $B_2$ are $0$ and $\delta_2=b_1-cb_2$. 
       The Lyapunov exponent is given by 
\be 
 \lambda=\frac{p_{21}}{p_{21}+p_{12}}\log |\delta_1|
+\frac{p_{12}}{p_{12}+p_{21}}\log |b_1-cb_2|+
\frac{p_{12}p_{21}}{p_{12}+p_{21}}\log| \frac{b_1}{b_1-cb_2}|.
\label{eq:lambda-ex3}
\ee
  We want to choose $b_1, b_2, c, \delta_1$ and  $p_{ij}$'s so that 
$|\delta_1|<1, |\delta_2|<1$ and $\lambda>0$. Since the first two terms in 
(\ref{eq:lambda-ex3}) are negative, we need to make the third term as large
as 
possible. Thus, $b_1/(b_1-cb_2)$ should be large. For example, if we choose

$b_1=100, c=10, b_2=9.99, \delta=0.1$.  Then in order $\lambda>0$ we
require 
\[ -p_{21} \log |\delta_1|+p_{12}\log 10 < 3p_{21} p_{12} \log 10.\]
This is satisfied if e.g. $\delta_1=0.1, p_{21}=p_{12}=0.8$.
 \hfill $\Box$

 If one subprocess is stable, the other is unstable, in most situations
there always exists a switching strategy to make the mixing process stable.
Consider the situation that  there exists a subordinate matrix norm such
that $\|B_1\|<1, \|B_2\|>1$. 
Then, $E\log \|A_1\|=\rho_1 \log \|B_1\|+\rho_2 \log \|B_2\|$ can be made
less than $0$ if $\rho_2$ is small enough. We call this strategy the
{\em preferred switching}, to denote the phenomenon that a mixture process with
less frequent unstable regime can  still be stable. 

Now we give an example
that two unstable vector processes 
can give rise to a stable mixing process.  

{\em Example 4.} Consider an extension of Example 2 to multidimensional
case when $B_I$'s commute. For example, let 
\[ B_1=\left( \begin{array}{cc} 2 & 0\\ 0 & \frac{1}{2} \end{array}
\right),
   B_2 =\left( \begin{array}{cc} \frac{1}{3} & 0 \\ 0 & \frac{3}{2}
\end{array} \right).\]
The two Lyapunov exponents associated with the switching between $B_1$ and
$B_2$ are given by 
\[\lambda_1=\rho_1\log 2- \rho_2\log 3, \lambda_2=-\rho_1\log2+\rho_2(\log
3-\log2).\]
We require that $\lambda_1<0$ and $\lambda_2<0$. Let $\rho=\rho_1$. 
This is true if and only if 
$$\frac{\log 3-\log 2}{\log 3}<\rho<\frac{\log 3}{\log 2+\log 3}.$$ 
 \hfill $\Box$

\subsection{\bf Mean shifting models}
\label{sec:shift}  
Consider the mean shifting model given by 
\be
X_n=M_n+A_n X_{n-1} +E_n
\label{eq:mean-shift}
\ee
where $A_n$ and $E_n$ as before and 
$M_n$ is the shifting mean, defined as $\mu_i$ when $I_n=i$ for 
$i=1,2,\ldots, r$ or $M_n=\sum_{i=1}^r \mu_i 1_{\{I_n=i\}}$.

The mean-shifting model can be regarded as a
more general case of SAR  when $E_n$ 
may be allowed to take nonzero mean as well, such as,
$\mu_i$ when $I_n=i$
for some $i$. An interesting case is when $A_n$ is a constant 
and only the mean or variance of $E_n$ varies among different 
regimes. 

Obviously $M_n$ is a stationary sequence if $I_n$ is. 
Using an expansion similar to (\ref{eq:expansion}) and 
Proposition~\ref{th:strict-stat}, it can be shown that the proper stationary 
solution of (\ref{eq:mean-shift}) is given by 
\be
W_n' =(M_n+\sum_{i=0}^{\infty} A_n A_{n-1} \cdots A_{n-i} M_{n-i-1})
    +(E_n+\sum_{i=0}^{\infty} A_n A_{n-1} \cdots A_{n-i} E_{n-i-1}).
\label{eq:limit-shift-mean}
\ee
That is, the stationary solution of (\ref{eq:mean-shift}) is given by 
the sum of two stationary processes 
\be
\bar{M}_n=M_n+\sum_{i=0}^{\infty} A_n A_{n-1} \cdots A_{n-i} M_{n-i-1}
\label{eq:limiting-mean}
\ee
and $W_n$ of (\ref{eq:limiting-proc}).
Note that (\ref{eq:limiting-mean}) is in general 
well-defined under the 
negative Lyapunov exponent assumption [cf. (\ref{eq:lyapunov-neg})] 
and  
\[\E \max(\log \|M_1\|,0)<\infty\]
(cf. Proof of Lemma~\ref{la:converg}).  
In particular, above condition is satisfied 
if $M_n$ takes on values from a finite set.

{\em Example 5.} Hamilton (1989)'s  model for  business cycle 
uses a fourth-order autoregression and mean-shifting model with two regimes.  
Writing in our  state space representation (\ref{eq:mean-shift}), 
this corresponds to $A_n$ taking on a fixed 
$A=\left( \begin{array}{llll} a_1 & a_2 & a_3 & a_4\\
                               1 & 0    & 0  & 0\\
                               0 & 1  & 0 & 0\\
                               0 & 0 & 1&0 \end{array} \right)$
and $M_n$ taking on $(\mu_i, 0, 0, 0)^T$ depending on $I_n=i, i=1,2$. 
By our theory, this model has a stationary  and stable solution as long as 
$A$ is stable. In particular, the empirical model  of 
Krolzig (1997, Sec.11.3.1) 
for German business cycle with $a_1=0.2932, a_2=0.1055, 
a_3=0.0026, a_4=0.3812$ clearly has a stationary solution  because  
$a_1+a_2+a_3+a_4<1$ and $a_i$'s are positive.

{\em Example 6.} We discuss another mean-shifting model 
which is a simplified version of  Lu and Berliner (1997)'s 
model for a riverflow time series $y_n$. 
 Their model 
consists of  mixture of   AR(1), ARX(1), and  
AR(1)  models with different means at each of  
the three regimes, corresponding to normal ($0$), 
rising (1), or falling (2) of 
the riverflow, where in the rising regime the past rainfall $x_{n-1}$   
series is included linearly.    We assume here that the 
regime switching process is independent of both $\{x_n\}, \{y_n\}$ 
and follows a Markov chain.  This model can be easily embedded in our 
formulation (\ref{eq:mean-shift}) with $p=1$ and $M_n$ taking on 
fixed values except in the rising regime when $M_n=\mu_1+a x_{n-1}$. 
Extending slightly the argument used in this section, if the 
rainfall series $\{x_n\}$  is {\em stationary} and the regime 
switching process is {\em ergodic}, the riverflow 
series $\{y_n\}$ is stationary if the AR(1) processes 
are either stationary or nonstationary of the  random walk type 
(cf. Example 2).

\section{Existence of moments}
\label{sec:second-order}
Existence of moments is often assumed in time series analysis, 
notably for the second-order theory (cf. Brockwell and Davis, 1991). 
For a general  stochastic difference  equation, Karlsen (1990) 
gives some  general conditions for checking the 
existence of finite moments. He also gives some examples 
where more explicit results can be derived. 
In this section,  by exploiting the Markovain structure in the 
hidden state process, we derive directly some explicit conditions for 
existence of   second-order moment  of SVAR models and the 
related autocorrelation property. 
 
We make the following assumption.\\
({\bf A}) $\lim_{n\arrow \infty}^{} \E [\|A_n \ldots, A_1\| | I_0=i]= 0$ 
for any $i=1,\ldots, R$.

By the ergodicity of $\{I_n\}$, one can easily show that ({\bf A}) is equivalent 
to 
the condition that \\
({\bf A}') $\lim_{n\arrow \infty}^{} \E \|A_n \ldots, A_1\|=0$.

Consider the property of the quantity defined by 
\[\Phi_{ni}(I_i)=E[\|A_n\ldots, A_{i+1}\| | I_i] \]
for any $n, i<n$. Then, since $\{A_n\}$ is an induced 
matrix-valued FMC defined in terms  of $I_n$. It shares the usual 
Markov property, and in particular 
$\Phi_{ni}(I_i)$ is independent of $i$ and depends only on $n-i$.
If we write 
\[ \Phi_{\ell}(I_0)=\E[\|A_{\ell}\ldots, A_{1}\| | I_0]\]
then
\[ \Phi_{ni}(I_i)=\Phi_{n-i}(I_i).\]
We use $\Phi_{\ell}$ or $\Phi_{ni}$ to denote their 
unconditional analogues. 

We have the following proposition on $\Phi_{\ell}(I_0)$. 
\begin{proposition}\label{prop:phi-func}
\[\mbox{ $\Phi_{n}(I_0)\arrow 0$   
if and only if $\Phi_{n}(I_0)$ tends to $0$ 
geometrically}.\] 
\end{proposition}
Proof: Since $\Phi_{n}(I_0)\arrow 0$ uniformly over $I_0$. Then, there exist 
an integer $\ell$ and constant   $\gamma<1$ such that 
$\Phi_{\ell}(i)\leq \gamma$ for all $i$. 
\[\Phi_{2\ell}(I_0) \leq 
\E[\Phi_{\ell}(I_{\ell}) \|A_{\ell}\ldots A_1\| | I_0]\]
\[\leq \gamma \Phi_{\ell}(I_0)\leq \gamma^2\]
There exists a constant $C$ such that 
$\Phi_n\leq C \gamma^{[n/\ell]}$ for any $n$. 
That is, $\Phi_{n}$ tends to $0$ at a geometric rate. 
The sufficient part of the proof is easy to establish.
\hfill $\Box$

\begin{theorem}\label{th:m-swtch-1}
The SVAR process  has a stationary solution whose second-order moment 
exists if (A) is satisfied.
\end{theorem}
Proof:
Consider the expansion for SVAR in (\ref{eq:stochastic-eq}):
\[
 X_n =A_n\ldots A_1 X_0 +A_n \ldots A_2 E_1
          +\cdots
         +A_n E_{n-1}+E_n.
\]
Then,
\beqa 
 E\|X_n\| & \leq &\E\|A_n\ldots A_1 \| \cdot \|X_0\| 
         +\E\|A_n \ldots A_2\|\cdot \|E_1\|\\
         &&+ \cdots
         +\E\|A_n\|\cdot \|E_{n-1}\|
          + \E\|E_n\|\\
         &=& \E \{\E[\|A_n\ldots A_1 \| | I_0]\cdot \|X_0\|\}
          +\E \{\E[\|A_n \ldots A_2\| I_1]\} \cdot 
          \|\Sigma_{I_1} \e_{nI_1}\|\}\\
         && + \cdots+
\E \{\E[\|A_n\| | I_{n-1}]\cdot \|\Sigma_{I_{n-1}}\e_{n-1 I_{n-1}}\|\}
          + \E\|E_n\|\\
        &\leq&\max \Phi_n(i) \E \|X_0\|
         + \max\Phi_{n-1}(i) 
          \cdot \E\|E_{1}\|\\
        &&+\cdots +\max \Phi_1(i) \cdot \E \|E_{n-1}\|+\E\|E_n\|
\eeqa
which is convergent if $\E\|X_0\|<\infty$, 
by Proposition~\ref{prop:phi-func} and 
ergodicity of $\{I_n\}$. 
Here assumptions on $\{E_n\}$ and 
 independence of $\{I_n\}$ and $\{\e_{ni}\}$ are  used.
\hfill $\Box$

Note that, by the concave nature of $\log X$, the Jensen's Inequality 
implies that the strict inequality 
\be 
\E \log \|A_n\ldots A_1\|< \log \E\|A_n\ldots A_1\|
\label{eq:jensen-log}
\ee
holds.

We denote  ${\limsup}_{n\arrow \infty} (1/n) 
\log \E \|A_n\ldots A_1\|$ by $\log(\gamma)$. 
Condition ({\bf A}) is equivalent to $\gamma<1$. By (\ref{eq:jensen-log}), 
this further implies that 
\be
  \lambda< \log \gamma<0.
\ee
This indicates that condition ({\bf A}) or 
({\bf A}') is stronger than 
negativity of the largest Lyapunov exponent $\lambda$, 
a potentially  general condition for strict stationarity. 
However, the latter  does not  
even ensure existence of second-order moment,  
see Bougerol and Picard for an example in the case of 
an GARCH process. 

Using the fact that  
\[X_{n+m} =A_{n+m} \ldots A_{n+1} X_n + A_{n+m}\ldots A_{n+2} E_{n+1} 
+ \ldots A_{n+m} E_{n+m-1} +E_{n+m}\]
for any integers $m$ and $n$,  
we have 
\beqa
 |\E X_n^T X_{n+m}| & =&|\E X_n^T A_{n+m} \ldots A_{n+1} X_n|\\
   &  \leq & \E|<X_n, A_{n+m} \ldots A_{n+1} X_n>| \\
   & \leq & \E \|A_{n+m} \ldots A_{n+1}\| \cdot  \|X_n\|^2.
\eeqa
That is, 
\be
 |\E X_n^T X_{n+m}| \leq \Phi_m \E \|X_1\|^2
\label{eq:aucovariance}
\ee
where we use the property that $\{X_n\}$ is {\em causal} and stationary, and 
$\{A_{n+i}\}$ is stationary.  
Thus,  the autocovariance matrix  at lag $m$ of the vector 
time series $\{X_n\}$ 
decays at a geometric rate, and  is bounded by $\gamma^m$.
 
\section{Switching ARMA models}
\label{sec:sarma}
We note some extensions of the switching autoregressive 
models. First, a switching moving average process of order
 $q$ (SMA(q)) can be defined as 
\be
X_n=E_n+ C_{1n} E_{n-1}+C_{2n} E_{n-2}+\ldots+C_{qn} E_{n-q}
\label{eq:sma}
\ee 
where $\{E_n\}$ is defined as before, and 
$E_{n-j}=\sum_{i=1}^r \Sigma_i \e_{(n-j) i} 1_{\{I_n=i\}}$ for 
$j=1,2,\ldots, q$. 
The coefficient matrices $C_{jn}$ will take on  member of a set of 
$r$ matrices depending on the value of $I_n$ for 
each $j$ between $1$ and $q$. 

We also assume that $\{ (\e_{n1},\ldots, \e_{nr})^T \}$ is 
stationary as before. 
If  $\{I_n\}$ is  stationary,   it follows that $\{X_n\}$ is stationary 
since it is a moving average function of stationary processes. 
On the other hand, if $\{I_n\}$ is ergodic, for arbitrary starting regime, 
$\{I_n\}$ eventually converges to stationarity and 
thus $\{X_n\}$ is asymptotically 
stationary. 

Similar to ARMA models, one can define switching ARMA (SARMA) models in which 
the coefficient matrices in both AR part and MA part take on 
different values depending on the {\em current} regime.  The 
stationarity condition for SVAR(1) models is also sufficient  
for SARMA(1,q) models. Since an AR(p) process can be represented 
as a vector AR(1) process, our theory applies to any switching 
ARMA(p,q) process. 

Other extension is also possible. In particular, 
the transition probabilities  of switching may be allowed to 
depend on past values of the process, or past values of 
another process.  
This  interesting  class of nonlinear 
time series models is closely related to  some traditional 
state dependent nonlinear time series models (cf. Tong 1990). 
Not surprisingly, there are increasing interest in 
applying them in some real modelling situations  
such as  security time series and  high-frequency data.  
It is our hope that the present work may shed light 
on these more complex models.

\end{document}